\documentclass[review]{elsarticle}

\usepackage{lineno,hyperref, amsmath,amsthm,amsfonts, enumitem}
\modulolinenumbers[5]

\journal{Journal of \LaTeX\ Templates}

\theoremstyle{definition}

\begin{document}

\begin{frontmatter}

\title{The comb-like representations of cellular ordinal balleans}

\author{I.V. Protasov}
\author{K.D. Protasova}

\begin{abstract}

Given two ordinal $\lambda$ and $\gamma$, let $f:[0,\lambda) \rightarrow [0,\gamma)$
be a function such that, for each  $\alpha<\gamma$, $\sup\{f(t): t\in[0, \alpha]\}<\gamma.$
We define a mapping $d_{f}: [0,\lambda)\times [0,\lambda) \longrightarrow [0,\gamma)$ by the rule: if $x<y$ then $d_{f}(x,y)= d_{f}(y,x)= \sup\{f(t): t\in(x,y]\}$, $d(x,x)=0$. The pair $([0,\lambda), d_{f})$ is called a $\gamma-$comb defined by $f$. We show that each cellular ordinal  ballean can be represented as a $\gamma-$comb.  In {\it General Asymptology}, cellular ordinal  balleans play a part of ultrametric spaces.



\end{abstract}

\end{frontmatter}


\section{Introduction}

In   \cite{b3}, a function $f: [0,1]\rightarrow [0,\infty)$ is called a {\it comb} if, for every $\varepsilon>0$, the set $\{t\in [0,1]: f(t)\geq \varepsilon\}$ is finite. Each comb $f$ defines a pseudo-metric $d_{f}$ on the set $I_{f}=\{t\in[0,1]: f(t)=0 \}$ by the rule: if $x<y$ then $$d_{f}(x,y)=\max\{f(t): t\in(x,y) \}, \ \  d_{f}(y,x)=d_{f}(x,y), \ \  d(x,x)=0.$$

After some reduced completion of $(I_{f}, d_{f})$, the authors get a compact ultrametric space and show that each compact ultrametric space with no isolated points can be obtained in this way.

In this note, we modify the basic construction from \cite{b3} to get the comb-like representations of cellular ordinal balleans which, in {\it General Asymptology} \cite{b7}, play a part of ultrametric spaces.

\section{Balleans}

Following \cite{b5}, \cite{b7}, we say that a {\em ball structure} is a triple $\mathcal{B}=(X,P,B)$, where $X$, $P$ are non-empty sets, and for all $x\in X$ and $\alpha\in P$, $B(x,\alpha)$ is a subset of $X$ which is called a {\em ball of radius} $\alpha$ around $x$. It is supposed that $x\in B(x, \alpha)$ for all $x\in X$, $\alpha\in P$.
The set $X$ is called the {\it support} of $\mathcal{B}$, $P$  is called the {\it set of radii}.

Given any $x\in X$, $A\subseteq X$, $\alpha\in P$, we set
$$B^*(x,\alpha)=\{y\in X:x\in B(y,\alpha)\},\ B(A,\alpha)=\bigcup_{a\in A}B(a,\alpha),\ B^*(A,\alpha)=\bigcup_{a\in A}B^*(a,\alpha).$$

A ball structure $\mathcal{B}=(X,P,B)$ is called a {\it ballean} if

\begin{itemize}
\item{} for any $\alpha,\beta\in P$, there exist $\alpha',\beta'\in P$ such that, for every $x\in X$,
$$B(x,\alpha)\subseteq B^*(x,\alpha'),\ B^*(x,\beta)\subseteq B(x,\beta');$$

\item{} for any $\alpha,\beta\in P$, there exists $\gamma\in P$ such that, for every $x\in X$,
$$B(B(x,\alpha),\beta)\subseteq B(x,\gamma);$$

\item{} for any $x,y\in X$, there exists $\alpha\in P$ such that $y\in B(x, \alpha)$.
\end{itemize}

We note that a ballean can be considered as an asymptotic counterpart of a uniform space, and could be defined \cite{b8} in terms of the entourages of the diagonal $\Delta_X = \{(x,x): x\in X\}$ in $X\times X$. In this case a ballean is called a {\em coarse structure}.

For categorical look at the balleans and coarse structures as "two faces of the same coin" see \cite{b2}.

Let $\mathcal{B}=(X,P,B)$, $\mathcal{B'}=(X',P',B')$ be balleans.
A mapping $f:X\to X'$ is called a $\prec$-\emph{mapping} if, for every $\alpha\in P$, there exists $\alpha'\in P'$ such that, for every $x\in X$, $f(B(x,\alpha))\subseteq B'(f(x),\alpha')$.

A bijection $f:X\rightarrow X'$ is called an {\it asymorphism} between $\mathcal{B}$ and $\mathcal{B}'$ if $f$ and $f^{-1}$ are $\prec$-mappings. In this case $\mathcal{B}$ and $\mathcal{B}'$ are called {\it asymorphic}.

 Given a ballean  $\mathcal{B} = (X,P,B)$, we  define a preodering $<$ on $P$ by the rule: $\alpha<\beta$ if and only if $B(x,\alpha)\subseteq B(x,\beta)$ for each $x\in X$. A subset $P^{\prime}$ of $P$ is called  {\it cofinal} if, for every $\alpha\in P$,  there exists $\alpha^{\prime}\in P^{\prime}$ such that $\alpha<\alpha'$. A ballean  $\mathcal{B}$ is called {\it ordinal} if there exists a cofinal well-ordered (by $<$) subset $P^{\prime}$ of $P$.

For  a ballean  $\mathcal{B} = (X,P,B)$, $x,y\in X$ and $\alpha\in P$, we  say that $x$ and $y$ are $\alpha${\it -path connected} if there exists a finite sequence $x_{0}, \ldots, x_{n}, \  x_{0}=x, x_{n}=y$ such that $x_{i+1}\in B(x_{i}, \alpha)$ for each $i\in \{0,\ldots,n-1\}$. For  any $x\in X$ and $\alpha\in P$, we set
$$B^{\diamond} (x,\alpha)=\{y\in X: x,y \  are \  \  \alpha-path \  connected \},$$
and say that the ballean $\mathcal{B}^{\diamond} = (X,P,B^{\diamond})$ is a {\it cellularization} of $\mathcal{B}$. A ballean $\mathcal{B}$ is called  {\it cellular} if the identity $id: X\rightarrow X$ is an asymorphism between $\mathcal{B}$ and $\mathcal{B}^{\diamond}$.

Each metric space $(X,d)$ defines a metric ballean  $\mathcal{B}(X,d) = (X,\mathbb{R}^{+},B_{\alpha})$, where $B_{d} (x,r)=\{y\in X: d(x,y)\leq r\}. $ Clearly, $\mathcal{B}(X,d)$ is ordinal and, if $d$ is an ultrametric then $\mathcal{B}(X,d) $ is cellular.

For examples, decompositions and classification of cellular ordinal balleans see \cite{b1},\cite{b2}, \cite{b4}, \cite{b6}.

\section{Representations}

For ordinals $\alpha, \beta$,  $  \alpha<\beta$, we use the standard notations

$ [\alpha, \beta]=\{t: \alpha\leq t\leq \beta\},   [\alpha,\beta)=\{t: \alpha\leq t<\beta\},$
$( \alpha ,\beta]=\{t: \alpha< t\leq\beta\} .$

Let $X$ be a set and $\gamma$ be an ordinal. We say that a mapping $d: X\times X\longrightarrow  [0,\gamma)$ is a $\gamma$-{\it ultrametric} if $d(x,x)=0$, $d(x,y)=d(y,x)$ and  $$d(x,y)\leq \max \{d(x,z), d(z,y)\}.$$

Clearly, each ultrametric space with integer valued metric is an $\omega$-ultrametric space. By [7, Theorem 3.1.1], every cellular metrizable ballean is asymorphic to some $\omega$-ultrametric space.

Given two $\gamma$-ultrametric spaces $(X,d)$, $(X^{\prime},d^{\prime})$, a bijection $h: X \longrightarrow X^{\prime}$ is called an {\it isometry} if, for any $x,y\in X$, we have $d(x,y)=d^{\prime}(h(x),h(y))$.

Now let $\lambda,\gamma$ be ordinal and $f: [0, \lambda)\rightarrow [0, \gamma)$ be a function such that, for each $\alpha<\lambda$, $\sup\{f(t): t\in[0,\alpha]\}<\gamma$. We define a mapping $d_{f}: [0,\lambda)\times [0,\lambda)\rightarrow [0,\gamma)$ by the rule: if $x<y$  then $$d_{f}(x,y)=d_{f}(y,x)=\sup\{f(t): t\in (x,y]\},     d(x,x)=0,$$
and say that $([0,\lambda), d_{f})$ is a $\gamma$-{\it comb} determined by $f$.
Evidently, each $\gamma$-comb is a $\gamma$-ultrametric space.

{\bf Theorem. }
{\it Every $\gamma$-ultrametric space $(X, d)$ is isometric to some $\gamma$-comb $([0,\lambda), d_{f}).$}

\begin{proof}
We proceed on induction by $\gamma$. For $\gamma=1$, we just enumerate $X$  as $[0,\lambda)$ and take $f\equiv 0$. Assume that we have proved the statement for all ordinals less than $\gamma$ and consider two cases.

{\it Case 1.}
Let $\gamma$ is not a limit ordinal, so $\gamma=\gamma^{\prime}+1$. We partition $X=\bigcup\{X_{\delta}: \delta\in[0,\nu)\}$ into classes of the equivalence $\sim$ defined by $x\sim y$ if and only if $d(x,y)< \gamma^{\prime}$. If $\delta<\delta^{\prime}<\nu$ and $x\in X_{\delta}$, $y\in X_{\delta^{\prime}}$ then $d(x,y)= \gamma^{\prime}$.

 By the inductive hypothesis, each $X_{\delta}$ is isometric to some $\gamma^{\prime}$-comb $([0,\lambda_{\delta}),d_{f_{\delta}})$. We replace  inductively each $\delta\in[0,\nu)$ with consecutive intervals $\{[l_{\delta}, l_{\delta}+ \lambda_{\delta}): \delta\in[0,\nu)\}$, $l_{0}=0$ and define a function $f: [0,\lambda)\rightarrow [0,\gamma), [0,\lambda)= \bigcup\{[l_{\delta}, l_{\delta}+ \lambda_{\delta}): \delta\in[0,\nu)\}$  as follows. We put $f=f_{0}$ on $[0,\lambda_{0})$. If $\delta>0$ then we put $f(l_{\delta})= \gamma^{\prime}$ and $f(l_{\delta}+x)= f_{\delta}(x)$ for $x\in(0, \lambda_{\delta})$.

After  $\mid\nu\mid$ steps, we  get the desired $\gamma$-comb $([0,\lambda), d_{f})$.

{\it Case 2.} $\gamma$ is a limit ordinal. We fix some $x_{0} \in X$ and, for each
 $\delta <\gamma$, denote $X_{\delta} =\{x\in X: d(x_{0}, x)<\delta\}$. By the inductive hypothesis, there is an isometry $h_{\delta} : X_{\delta}\longrightarrow ([0, \lambda_{\delta}), d_{f_{\delta}}).$ Moreover, in view of Case 1, $f_{\delta + 1}$ and $h_{\delta + 1}$ can be chosen as the extensions of $f_{\delta}$ and $h_{\delta}$. Hence, we can use induction by $\delta$ to get the desired $\gamma$-comb and isometry.
\end{proof}

 Every $\gamma$-ultrametric  space $(X, d)$ can be considered as the ballean
 $(X, [0, \gamma), B_{d})$, where $B_{d}(x,\alpha)=\{y\in X: d(x,y)\leq \alpha\}. $

On the other hand, let $(X, P, B)$ be a cellular ordinal ballean. We may suppose that $P=[0,\gamma)$ and $B(x,\alpha)=B^{\diamond} (x,\alpha)$ for all $x\in X$, $\alpha\in [0,\gamma)$. We define a $\gamma$-ultrametric $d$ on $X$ by $d(x,y)= \min \{\alpha\in [0,\gamma): y\in B(x,\alpha)\}$. Then $(X,P,B)$ is asymorphic to $(X,d)$.

{\bf Corollary}
{\it Every cellular ordinal ballean is asymorphic to some $\gamma$-comb.
}

Taras Shevchenko National University of Kyiv, Department of Cybernetics, Volodymyrska 64, 01033, Kyiv Ukraine

E-mail address: i.v.protasov@gmail.com

E-mail address: k.d.ushakova@gmail.com

\bibliography{mybibfile}

\end{document}